\documentclass[11pt]{amsart}

\usepackage{amsmath,amsfonts,amssymb}
\usepackage{array,delarray}
\usepackage[dvips]{graphicx}

\newcommand{\al}{\ensuremath{\alpha}}
\newcommand{\be}{\ensuremath{\beta}}

\newcommand{\codim}{\mathrm{codim}}

\newcommand{\rad}{\mathrm{rad}}
\newcommand{\lcm}{\mathrm{lcm}}

\newcommand{\dms}{\Delta_{M^{\ast}}}
\newcommand{\dm}{\Delta_{M}}
\newcommand{\nn}{\mathbb{N}}

\newtheorem{theorem}{Theorem}[section]   

\newtheorem{corollary}[theorem]{Corollary}  

\newtheorem{lemma}[theorem]{Lemma}         

\newtheorem{proposition}[theorem]{Proposition}  

\theoremstyle{definition}
\newtheorem{definition}[theorem]{Definition}   

\theoremstyle{remark}

\newtheorem{example}[theorem]{Example}        

\numberwithin{equation}{section}     

\begin{document}

\title[Generic Cohen-Macaulay Monomial Ideals]        
{Generic Cohen-Macaulay Monomial Ideals}

\author[Jarrah and Laubenbacher]{Abdul Salam Jarrah and Reinhard Laubenbacher}

\address[Abdul Salam Jarrah]
{Department of Mathematics\\
East Tennessee State University\\
Johnson City, TN 37614} \email{jarrah@mail.etsu.edu}

\address[Reinhard Laubenbacher]
{
Virginia Bioinformatics Institute\\
1880 Pratt Drive\\
Blacksburg, VA 24061} \email{reinhard@vbi.vt.edu}

\thanks{}

\date{\today}

\begin{abstract}
Given a simplicial complex, it is easy to construct a generic
deformation of its  Stanley-Reisner ideal. The main question under
investigation in this paper is how to characterize the simplicial
complexes such that their Stanley-Reisner ideals have
Cohen-Macaulay generic deformations.  Algorithms are presented to
construct such deformations for matroid complexes, shifted
complexes, and tree complexes.
\end{abstract}



\maketitle

\section{Introduction}
\hspace{\parindent} Let $S:= k[x_1,\dots,x_n]$ be a polynomial
ring in $n$ variables over a field $k$. A monomial ${\bf x^a} :=
x_1^{a_1}x_2^{a_2}\cdots x_n^{a_n}$ is called \textit{square-free}
if $a_i = 0$ or 1 for all $1 \leq i \leq n$. For monomials ${\bf
x^a}$ and ${\bf x^b}$, we say ${\bf x^a}$ \textit{strictly
divides} ${\bf x^b}$ if, for $1 \leq i \leq n$ and $a_i > 0$, we
have $a_i < b_i$. A \textit{monomial ideal} $M$ is an ideal
minimally generated by monomials, say $m_1,\dots,m_t$, (written as
$M = \langle m_1,\dots ,m_t \rangle$). The monomial ideal $M$ is
called \textit{square-free} if all $m_i$ are square free
monomials.

\begin{definition}\label{generic:def}
A monomial ideal $M = \langle m_1,\dots,m_t \rangle$ is called
\textit{generic} if for any two distinct minimal generators $m_i$
and $m_j$ with the same positive degree in some variable $x_h$,
there exists a third minimal generator $m_l$ such that $m_l$
strictly divides $\lcm(m_i,m_j)$.
\end{definition}
\begin{example}\label{ex1:generic}
The monomial ideal  $\langle x^2y^2,x^2z^2,xz \rangle \subset
k[x,y,z]$ is generic but $\langle xy^2,xz^2,xz \rangle \subset
k[x,y,z]$ is not generic.
\end{example}
This definition of generic monomial ideals has appeared in
\cite{MSY},  generalizing a definition in \cite{BPS}. See
\cite{Y1} for more information about generic monomial ideals.

\begin{definition}\label{simp-comp:def}
A simplicial complex $\Gamma$ on the vertex set
$[n]:=\{1,\dots,n\}$ is a non-empty collection of subsets of $[n]$
such that $\{i\} \in \Gamma$ for all $i \in [n]$, and whenever
$\tau \subset \sigma$ and $\sigma \in \Gamma$, then $\tau \in
\Gamma$. The elements of $\Gamma$ are called \textit{faces}, and
the \textit{dimension} of the face $\sigma  \in \Gamma$ is
$\dim(\sigma) := |\sigma|-1$. The dimension of $\Gamma$ is
$\dim(\Gamma) := \max\{\dim(\sigma) \, : \, \sigma \in \Gamma \}$.
Faces of dimension zero are called \textit{vertices} and faces
that are maximal under inclusion are called \textit{facets}. The
set of \textit{minimal non-faces} of $\Gamma$ is $\Sigma_\Gamma :=
\{ \sigma \subseteq [n] \, : \, \sigma \notin \Gamma \mbox{ but }
\sigma \setminus \{i\} \in \Gamma  \mbox{ for each } i \in \sigma
\}$.
\end{definition}

There is a one-to-one correspondence between square-free monomial
ideals in $S$ and simplicial complexes on $[n]$. Given a
simplicial complex $\Gamma$ on the vertex set $[n]$, the
\textit{Stanley-Reisner} ideal $I_{\Gamma}$ of $\Gamma$ is the
square-free monomial ideal $\langle x_{i_1}\cdots x_{i_s} \, : \,
\{i_1,\dots, i_s\} \in [n] \setminus \Gamma \rangle \subset S$.
Conversely, for a given non-zero square-free monomial ideal $M
\subset S$, there exists a simplicial complex $\Gamma$ on $[n]$
such that $I_{\Gamma} = M$.
\begin{definition}\label{deformation:def}
Let $I = \langle m_1,\dots,m_t \rangle$ be a square-free monomial
ideal. A \textit{deformation} of $I$ is a monomial ideal $M =
\langle m_1^*,\dots,m_t^*\rangle$ such that, for all $1\leq i \leq
t$, $m_i^* = x_{i_1}^{j_1}\cdots x_{i_s}^{j_s}$, where $m_i =
x_{i_1} \cdots x_{i_s}$ and $(j_1,\dots,j_s) \in (\nn\setminus
\{0\})^s$.
\end{definition}
Note that our definition of a deformation is different from the
ones given at \cite{BPS,MSY}. In particular, our deformations are
monomial ideals.

For a given simplicial complex $\Gamma$, one can easily find a
deformation $M$ of the Stanley-Reisner ideal $I_\Gamma$ such that
$M$ is a generic monomial ideal. In this paper, we characterize
simplicial complexes such that their corresponding Stanley-Reisner
ideals have deformations that are generic and Cohen-Macaulay. To
be precise, we provide algorithms to construct generic
Cohen-Macaulay deformations from the Stanley-Reisner ideals of
matroids complexes, shifted complexes, and tree complexes. A
version of this question appeared in \cite{MSY} where generic
Cohen-Macaulay monomial ideals have been studied extensively.
\begin{definition}\label{scarf:def}
Let $M = \langle m_1,\dots,m_t \rangle$ be a monomial ideal in
$S$. For $\sigma \subset [t]$, let $m_{\sigma} := \lcm(m_i \, : \,
i \in \sigma)$. The \textit{Scarf complex} of $M$ is
\[
\dm := \{ \sigma \subseteq [t] \, : \, m_\sigma \neq m_\tau \mbox{
for all } \tau \subseteq [t] \mbox{ and } \tau \neq \sigma \}.
\]
For $D$ large enough (larger than any exponent of any variable in
any minimal generator of $M$), let
\[
M^{\ast} := M + \langle x_1^D, \dots, x_n^D \rangle.
\]
The \textit{extended Scarf complex} of $M$ is the Scarf complex
$\dms$ of $M^{\ast}$.
\end{definition}
It is easy to see that the induced subcomplex of $\dms$ on the set
of generators of $M$ is $\dm$.  On the other hand, Miller
\textit{et. al} \cite{MSY} proved that if $M$ is generic, then the
induced subcomplex of  $\dms$ on  the set $\{x_1^D, \dots,
x_n^D\}$ is the simplicial complex $V(M)$, where $I_{V(M)} =
\rad(M)$.
\begin{example}\label{gen2:ex1}
Let $\Gamma$ is the simplicial complex in Figure \ref{gen2:fig1},
and let $M = \langle x_1^2x_2^2x_3^2,x_1x_4, x_2^2x_4^2 \rangle$.
It is clear that $M$ is a generic deformation of $I_{\Gamma}$ (in
particular, $\rad(M) = I_{\Gamma}$ and hence $V(M) = \Gamma$). On
the other hand, let $M^{\ast} := M + \langle
x_1^3,x_2^3,x_3^3,x_4^3 \rangle$. Figure \ref{cm} gives $\dm$ and
$\dms$. The set of facets of $\dms$ is
\begin{align*}
\{\{x_1x_4,x_1^2x_2^2x_3^2,x_1^3,x_2^3\},\{x_1x_4&,x_1^2x_2^2x_3^2,x_1^3,x_3^3\},
 \{x_1x_4,x_1^2x_2^2x_3^2,x_2^3,x_3^3\}, \\
&\{x_1x_4,x_2^2x_4^2,x_2^3,x_3^3\},
\{x_1x_4,x_2^2x_4^2,x_3^3,x_4^3\}\}.
\end{align*}
\begin{figure}[!htp]
\centering
\includegraphics[totalheight=5cm]{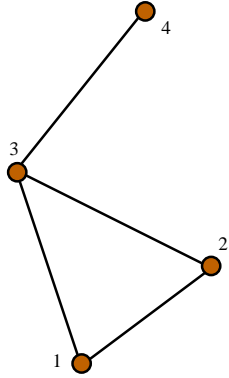}
\caption{The simplicial complex $\Gamma$ from Example
\ref{gen2:ex1}}. \label{gen2:fig1}
\end{figure}
\begin{figure}[!htp]
\centering
\includegraphics[totalheight=5cm]{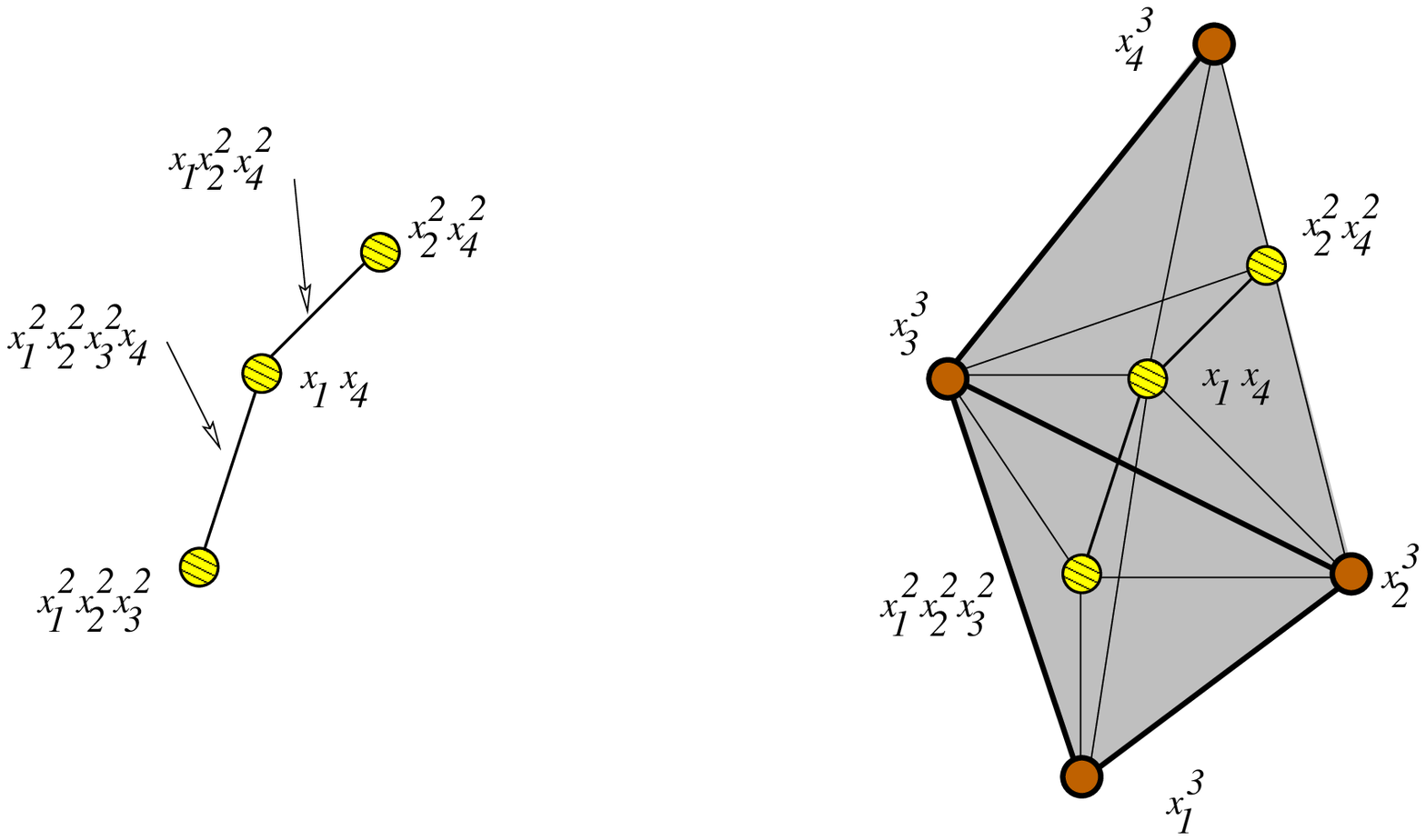}
\caption{$\dm$ and $\dms$  of $M$ in Example \ref{gen2:ex1}.}
\label{cm}
\end{figure}
\end{example}
The following theorem is proved as part of \cite[Theorem
1.7]{MSY}.
\begin{theorem}\label{primary:thm}
Let $M$ be a generic monomial ideal. For each facet $\sigma$ of
$\dms$, let
\[
M_{\sigma} := \langle x_s^{p_s} \, : \, p_s :=
\deg_{x_s}(m_{\sigma}) \mbox{ and } p_s < D \rangle.
\]
Then
\[
M = \bigcap_{ \sigma \mbox{ is a facet of } \dms} M_{\sigma}
\]
is a minimal irreducible decomposition of $M$.
\end{theorem}
\noindent \emph{Example} \ref{gen2:ex1} (\emph{cont.}).
\label{exam:generic} By Theorem \ref{primary:thm}, each facet of
$\dms$ corresponds to an irreducible component of $M$, and hence
$\dms$ gives a  minimal irreducible decomposition of $M$ which
implies a minimal primary decomposition of $M$:
\begin{align*}
M &= \langle x_3^2,x_4 \rangle \cap \langle x_2^2, x_4 \rangle
\cap \langle x_1^2,x_4 \rangle \cap \langle x_1,x_4^2 \rangle
\cap \langle x_1,x_2^2 \rangle \\
&= \langle x_3^2,x_4 \rangle \cap \langle x_2^2, x_4 \rangle \cap
\langle x_1^2,x_1x_4, x_4^2\rangle \cap \langle x_1,x_2^2 \rangle.
\end{align*}

\begin{theorem}[{\cite[Theorem 2.5]{MSY}}]\label{mainthm}
Let $M$ be a generic monomial ideal. Then $M$ has no embedded
associated primes if and only if  $M$ is Cohen-Macaulay. In  this
case, both $\dm$ and $V(M)$ are shellable.
\end{theorem}
For a given simplicial complex $\Gamma$, the theorem above implies
that if $\Gamma$ is not shellable, then there is no generic
Cohen-Macaulay monomial ideal $M$ such that $V(M)= \Gamma$. In
particular, there is no generic Cohen-Macaulay deformation of
$I_\Gamma$. The following theorem will play a major role
subsequently.
\begin{theorem} \label{thethm}
Let  $M \subset S$ be a  generic monomial ideal. Then $M$ is
Cohen-Macaulay if and only if $\dim(\dm) + \dim(V(M)) = n-2$.
\end{theorem}
\begin{proof}
Suppose $M$ is Cohen-Macaulay. By Theorem \ref{primary:thm}, each
facet of $\dms$ gives an irreducible component of $M$.  In
particular, for any facet $\sigma$ of $\dms$, we have
$$
|\sigma \cap \{m_1, \dots, m_t\}| = \codim(M) \;\;\;\; \mbox{and}
\;\;\;\; |\sigma \cap \{x_1^D, \dots, x_n^D\}| = \dim(R/M).
$$
Furthermore, both cardinalities are independent of the facet
$\sigma$, since $M$ has no embedded primes. But $\sigma \cap
\{m_1, \dots, m_t\}$ is a face of $\dm$ and any facet of
$\Delta_M$ is a restriction of a facet of $\dms$. Thus $\dim(\dm)
= \codim(M)-1$. Hence $\sigma \cap \{x_1^D, \dots, x_n^D\}$ is a
facet of $V(M)$ and therefore $ \dim(V(M)) = \dim(R/M)-1$. But
$\codim(M)+\dim(R/M) = n$. Thus, $\dim(\dm) + \dim(V(M)) = n-2$.

Conversely, suppose $\dim(\dm) + \dim(V(M)) = n-2$. Then  it is
enough to show that for any facet $\sigma$ of $\dms$, we have
$|\sigma \cap \{m_1, \dots, m_t\}| = n-d$, where $\dim(V(M)) =
d-1$, i.e., $\dim(\dm) = n-d-1$. This proves that all irreducible
components are of the same codimension, namely $n-d$, and hence
there are no inclusions among them. Thus the irreducible
decomposition is minimal and hence all the associated primes of
$M$ are of the same codimension. Thus $M$ has no embedded primes.
By  Theorem \ref{mainthm}, $M$ is then Cohen-Macaulay.

To this end, let $\sigma$  be a facet of  $\dms$. If $|\sigma \cap
\{m_1, \dots, m_t\}| > n-d$, then $\dim(\dm) > n-d-1$, a
contradiction. Now if $|\sigma \cap \{m_1, \dots, m_t\}| < n-d$,
then $|\sigma \cap \{x_1^D, \dots, x_n^D\}| > d$. Thus $\dim(V(M))
> d-1$, a contradiction. Therefore, $M$ is Cohen-Macaulay.
\end{proof}
\begin{example}\label{ex-6}
Let $\Gamma$ be the 1-dimensional simplicial complex (graph) in
Figure \ref{fig-6}.
\begin{figure}[!htp]
\centering

\includegraphics[totalheight=4cm]{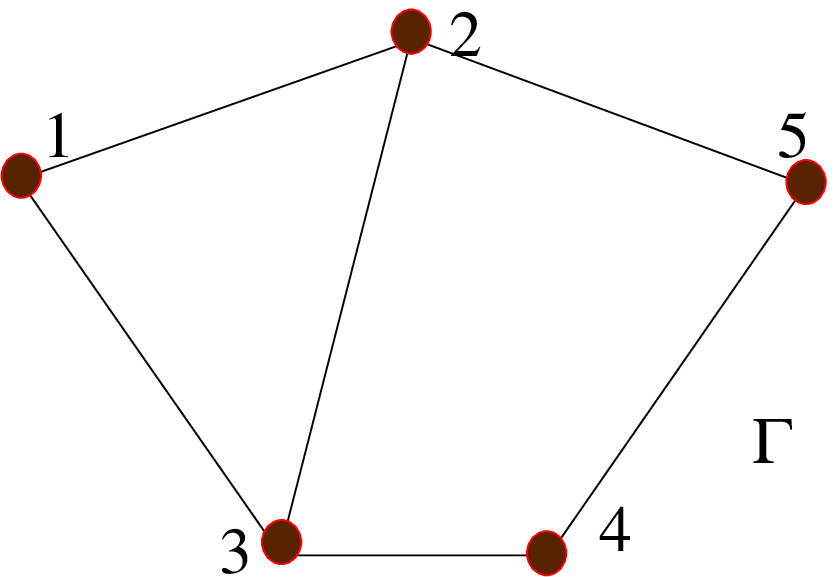}
\caption{The simplicial complex $\Gamma$ from Example \ref{ex-6}.}
\label{fig-6}
\end{figure}
It is clear that $\Gamma$ is shellable, since it is connected. The
Stanley-Reisner ideal of $\Gamma$ is
\[
I_{\Gamma}= \langle x_1x_2x_3, x_1x_4, x_1x_5, x_2x_4, x_3x_5
\rangle.
\]
Let $M_1 =\langle  x_1^3x_2^3x_3^3, x_2^3x_4^3\rangle$, $M_2
=\langle  x_1^2x_4^2, x_3^2x_5^2 \rangle$, $M_3 = \langle  x_1x_5
\rangle$, and let $M = M_1+M_2+M_3$. It is clear that $M$ is a
deformation of $I_{\Gamma}$. Moreover, $M$ is generic.  Namely,
for $1 \leq i \leq 2$, given any two distinct minimal generators
$m_1, m_2 \in M_i$, there exists a third generator $m_3$ in $M_j$
for some $j > i$ such that $m_3$ strictly divides $\lcm(m_1,m_2)$.
Thus $\dim(\dm) \leq 2 = 5-2-1=n-d-1$. But $\{ x_1x_5, x_1^2x_4^2,
x_2^3x_4^3\}$ is clearly a face of $\dm$. So $\dim(\dm)= 2$. Now,
by Theorem \ref{thethm}, $M$ is Cohen-Macaulay.
\end{example}
In the next sections, we show that the Stanley-Reisner ideals of
matroid, shifted, or tree complexes have generic Cohen-Macaulay
deformations.
\section{Matroid Complexes}
\hspace{\parindent} In this section, we prove that the
Stanley-Reisner ideals of matroid complexes have generic
Cohen-Macaulay deformations. There are many equivalent definitions
of matroids, see \cite{Ox}, or \cite[\S III.3]{St}. The following
definition uses the so-called \textit{circuits} axiom.
\begin{definition}\label{def1:matroid}
A simplicial complex $\Delta$ is a {\it matroid } if for any two
minimal non-faces (circuits)  $\al$ and $\be$ with an $i \in \al
\cap \be$, there exists a minimal non-face $\gamma$ such that
$\gamma \subseteq (\al \cup \be) \setminus \{i\}$.
\end{definition}
The following corollary follows directly from \cite[Proposition
3.1]{St}.
\begin{corollary}\label{cor1:matroid}
Let $\Gamma$ be a matroid complex on the vertex set $[n]$. For
every subset $W$ of $[n]$, the induced subcomplex $\Gamma_W := \{
\sigma \in \Gamma \, : \, \sigma \subseteq W \}$ of $\Gamma$ is a
matroid complex. In particular, $\Gamma_W$ is pure and shellable.
\end{corollary}
\begin{example}\label{ex1:matroid}
The simplicial complex $\Gamma$ in Figure \ref{fig1:matroid} is a
matroid. The Stanley-Reisner ideal of $\Gamma$ is $I_{\Gamma} =
\langle x_1x_3,x_1x_2x_4,x_2x_3x_4 \rangle$.
\begin{figure}[h]
  \centering
  \includegraphics[width=5cm]{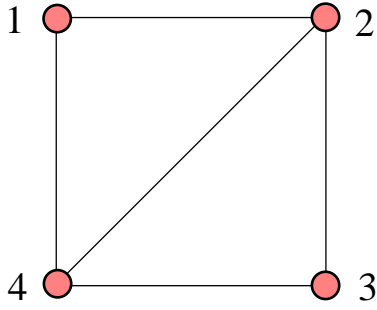}
  \caption{The matroid complex $\Gamma$ from  Example \ref{ex1:matroid}.}
  \label{fig1:matroid}
\end{figure}
\end{example}
The following theorem is the main result of this section.
\begin{theorem}\label{thm2:matroid}
Let $\Gamma$ be a matroid complex on the vertex set $[n]$ and
suppose $\dim(\Gamma) = d-1$. Then there exists a generic
Cohen-Macaulay deformation of $I_{\Gamma}$.
\end{theorem}
\begin{proof}
Pick a facet $\sigma$ of $\Gamma$. Without loss of generality,
$\sigma= \{1,\dots,d\}=[d]$. There exists a unique minimal
non-face $\al \subset [d+1]$: By Corollary \ref{cor1:matroid}, the
induced subcomplex $\Gamma_{[d+1]}$ is pure and $\sigma \in
\Gamma_{[d+1]}$, in particular, $\dim(\Gamma_{[d+1]}) = d-1$. But
$\Gamma_{[d+1]}$ is pure, so there exists a unique minimal
non-face $\al \subset [d+1]$. Let $m_\al = \prod_{j \in \al} x_j$.
Let  $M_1 := \langle m_\al \rangle$, a monomial ideal in $S$. By
induction, for $ i > 1$, let
\begin{align*}
\Sigma_{i} := \{ \al \subset [d+i] \, : \, \mbox{$\al$ is a
minimal non-face of $\Gamma$ and } d+i \in \al \},
\end{align*}
 be the set of minimal non-faces in $[d+i]$ but not
in $[d+i-1]$. Let
\begin{align*}
M_i := \langle \, \prod_{j \in \al} x_j^i \, : \, \al \in \Sigma_i
\,  \rangle.
\end{align*}
Since we have $n$ vertices, we will have all the minimal non-faces
at the inductive step $n-d$.

Let $M = M_1 + \cdots + M_{n-d}$. Then $M$ is generic and $V(M) =
\Gamma$. It remains to show that $M$ is Cohen-Macaulay. By Theorem
\ref{thethm}, it is enough to show that $\dim(\dm) = n - d-1$.

First, we show that, for all $i > 0$, $M_i \neq 0$: at each
inductive step $i$ , $ 2 \leq i \leq n-d$, $\Gamma_{[d+i]}$ is a
matroid complex on $d+i$ vertices and $\dim(\Gamma_{[d+i]}) =
d-1$. For if $M_i = 0$, there is no missing face with the vertex
$d+i$, in particular $\sigma \cup \{d+i\} \in \Gamma$. This is a
contradiction, since $\sigma$ is a facet of $\Gamma$. Furthermore,
$\dim(\dm) = n-d-1$: For each $1 \leq i \leq n-d$, let
\begin{align*}
F_i := \sigma \cup \{d+i\}.
\end{align*}
Then there exists a unique minimal non-face (circuit) $c_i$ of
$\Gamma$ such that $d+i \in c_i \subset F_i$. Let
\begin{align*}
m_i := \prod_{j \in c_i} x_j^i.
\end{align*}
Let $G =\{m_1,\dots,m_{n-d}\}$. We will show that $G$ is a face of
$\dm$ and hence $\dim(\dm) = n-d-1$. Notice that
\[
\lcm(G) := \lcm(m_1,\dots, m_{n-d})= f \cdot x_{d+1}x_{d+2}^2
\cdots x_{n}^{n-d},
\]
where $f$ is a monomial in the variables $x_1, \dots, x_d$.
Moreover, for any proper subset $H \subset G$, it is easy to see
that $\lcm (H) \neq \lcm (G)$.

Suppose that, for some $1 \leq t \leq n-d$, there exists a
generator
 $m \in M_t$ such that $m$ divides $\lcm(G)$,
say $m = x_{i_1}^tx_{i_2}^t\cdots x_{i_s}^tx_{d+t}^t$, where $1
\leq i_1 \leq \cdots \leq i_s < d+t$. Thus $i_s \leq d$. This
implies that $m$ corresponds to the unique minimal non-face $\al
=\{i_1, \dots, i_s, d+t\}$ which is $m_t$. Hence $G$ is a face of
$\dm$ and therefore $\dim(\dm) = n-d-1$.
\end{proof}
\begin{example}\label{ex2:matroid}
For the simplicial complexes $\Gamma$ in Figure
\ref{fig1:matroid}, $\sigma = \{1,2\}$ and $\Gamma_{[3]}$ has the
facets $\sigma$ and $\{2,3\}$. Thus $M_1 = \langle x_1x_3
\rangle$. It is clear that $\Gamma_{[4]} = \Gamma$. Thus we need
to compute only $M_2$. But $\Sigma_2 = \{\{1,2,4\}, \{2,3,4\}\}$.
Thus $M_2 = \langle x_1^2x_2^2x_4^2,,x_2^2x_3^2x_4^2 \rangle$. Let
$M = M_1 + M_2 = \langle  x_1x_3, x_1^2x_2^2x_4^2,,x_2^2x_3^2x_4^2
\rangle$. It is clear that $M$ is generic. Moreover,
$\dim(\Delta_M) = 1$, see Figure \ref{fig2:matroid}. Thus $M$ is
Cohen-Macaulay, by Theorem \ref{mainthm}.
\begin{figure}[h]
  \centering
  \includegraphics[width=5cm]{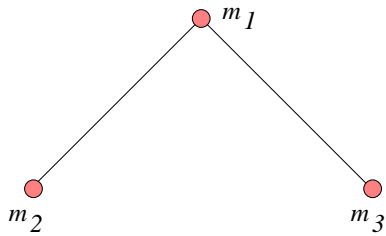}
  \caption{The Scarf complex $\dm$ of $M$ in Example \ref{ex2:matroid}.}
  \label{fig2:matroid}
\end{figure}
\end{example}
\section{Shifted Complexes}
\hspace{\parindent} In this section, we prove that shellable
shifted complexes have  generic Cohen-Macaulay deformations.
\begin{definition}\label{def1:shifted}
A simplicial complex $\Gamma$ on the vertex set $[n]$ is {\it
shifted } if, for all $\sigma \in \Gamma$, whenever $j \in
\sigma$,
 $i \in [n]$, and $i < j$,  we have
$(\sigma \setminus \{j\}) \cup \{i\} \in \Gamma$.
\end{definition}
\begin{example}\label{ex1:shifted}
Let $\Gamma$ be the 2-dimensional simplicial complex on the vertex
set $[6]$ with the set of facets $\{\{1,2,3\},
\{1,2,4\},\{1,3,4\},\{1,2,5\},\{1,2,6\}\}$, see Figure
\ref{fig1:shifted}. It is easy to see that $\Gamma$ is a shifted
complex.
\begin{figure}[h]
  \centering
 \includegraphics[width=4.0in]{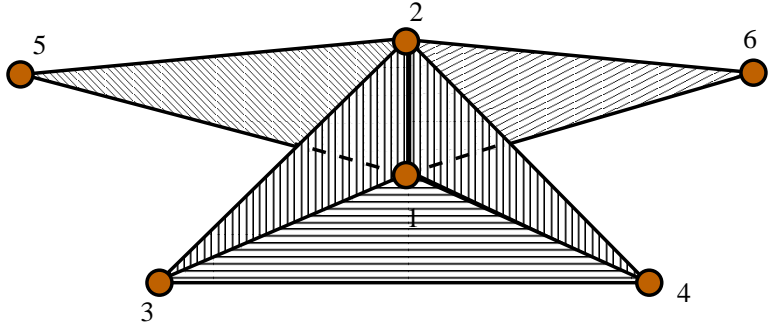}
  \caption{The shifted complex $\Gamma$ from Example \ref{ex1:shifted}.}
  \label{fig1:shifted}
\end{figure}
\end{example}
Notice that the labeling plays a major role in this definition. So
relabeling (permuting the vertices)
 a shifted complex might yield a non-shifted one, although
they are {\it combinatorially equivalent} (their face lattices are
isomorphic).

The following theorem \cite[Theorem 3]{BK1} gives simple
geometrical and combinatorial characterizations for
Cohen-Macaulay, and shellable shifted complexes.
\begin{theorem}\label{thm1:shifted}
Let $\Gamma$ be a $(d-1)$-dimensional shifted simplicial complex.
\begin{enumerate}
\item $\Gamma$ is homotopically equivalent to a wedge of spheres,
\item $\Gamma$ is shellable if and only if  $\Gamma$ is
Cohen-Macaulay if and only if  $\Gamma$ is pure if and only if for
every $\sigma \in \Gamma$ and $|\sigma| = t$, $\sigma \cup [d-t]
\in \Gamma$.
\end{enumerate}
\end{theorem}
Let $\Gamma$ be a shellable $(d-1)$-dimensional shifted simplicial
complex on the vertex set $[n]$. We give a partition for the set
of minimal non-faces $\Sigma_{\Gamma}$. Recall that
\[
\Sigma_\Gamma :=  \{ \sigma \subseteq [n] \, : \, \sigma \notin
\Gamma \mbox{ but } \sigma \setminus \{i\} \in \Gamma  \mbox{ for
each } i \in \sigma \}.
\]
Now, for $1 \leq j \leq n-d$, let
\[
\Sigma_j := \{ \sigma \in \Sigma_{\Gamma} \, : \, \sigma = A \cup
\{d-|A|+j\} \mbox{ and } A \subseteq \overline{[d-|A|+j]}\}.
\]
In  the next proposition, we present some facts that will be very
useful in proving the main result of this section. But before that
we give an example.

\vspace{.2cm} \noindent {\it Example} \ref{ex1:shifted} ({\it
cont.}){\bf .} It is clear that
\[
\Sigma_\Gamma = \{\{2,3,4\}, \{3,5\}, \{3,6\}, \{4,5\}, \{4,6\},
\{5,6\}\}.
\]
Moreover, it is easy to see that
\begin{align*}
\Sigma_1 &=  \{\{2,3,4\}, \{3,5\}, \{3,6\}\}. \\
\Sigma_2 &=  \{\{4,5\}, \{4,6\}\}. \\
\Sigma_3 &=  \{\{5,6\}\}.
\end{align*}
Hence, $\Sigma_\Gamma = \Sigma_1 \cup \Sigma_2 \cup \Sigma_3$.
\begin{proposition}\label{thm2:shifted}
Let $\Gamma$, $\Sigma_{\Gamma}$, and $\Sigma_j$ be as above.
\begin{enumerate}
\item For $1 \leq j \leq n-d$,  we have $\Sigma_j \neq \emptyset$,
\item For $1 \leq i < j \leq n-d$, $\Sigma_i \cap \Sigma_j =
\emptyset$, and \item  $\Sigma_{\Gamma} = \cup_{j=1}^{n-d}
\Sigma_j$.
\end{enumerate}
\end{proposition}
\begin{proof}
1. Let  $1 \leq j \leq n-d$ and let $H:= \{d+j,d+j-1,\dots, j\}$.
Since $|H| = d+1$ and $\dim(\Gamma)=d-1$, for some  $j \leq t \leq
d+j$, we have that $\{d+j, d+j-1, \dots, t\}$ is a minimal
non-face of $\Gamma$. Let $A := \{d+j, d+j-1, \dots,t-1\}$. It is
clear that $|A| = d-t+j$ and hence $t = d-|A|+j$. Therefore, $A
\cup \{t\} \in \Sigma_j$ which implies that $\Sigma_j \neq
\emptyset$.

2. Straightforward.

3. It is clear that, for any $1 \leq j \leq n-d$, we have
$\Sigma_j \subseteq \Sigma_{\Gamma}$. Let $\sigma \in
\Sigma_{\Gamma}$, say $\sigma = \{i_1, \dots, i_{t+1} \}$ where $1
\leq i_{t+1} < \cdots < i_1 \leq n$. Let $A := \{i_1,\dots,
i_t\}$. Since $\sigma$ is a minimal non-face, we have $A \in
\Gamma$. Since $1 \leq i_{t+1} \leq n-t$, there exists $ 1 \leq
j_0 \leq n-d$ such that $i_{t+1} = d-t+j_0$. Notice that $|A|= t$.
Thus $A \subseteq \overline{[d-|A|+j_0]}$ and $\sigma = A \cup
\{d-|A|+j_0\}$. Therefore $\sigma \in \Sigma_{j_0}$. Hence
$\Sigma_{\Gamma} = \cup_{j=1}^{n-d} \Sigma_j$.
\end{proof}
\begin{proposition}\label{thm3:shifted}
Let $\Gamma$, $\Sigma_j$ be as above. Let $\sigma, \gamma \in
\Sigma_j$ and $\sigma \neq \gamma$. There exists $ s > j$ and $
\tau \in \Sigma_s$ such that $\tau \subset \sigma \cup \gamma$.
\end{proposition}
\begin{proof}
Since $\sigma, \gamma \in \Sigma_j$, there exists $A, B \in
\Gamma$ such that $\sigma = A \cup \{d-|A|+j\}$ and $ \gamma = B
\cup \{d-|B|+j\}$ where $A \subseteq \overline{[d-|A|+j]}$ and $B
\subseteq \overline{[d-|B|+j]}$. Moreover, $A \neq B$, since
$\sigma \neq \gamma$. Without loss of generality, assume $|A| \leq
|B|$.

\noindent {\bf Case 1}. Suppose $b \leq d-|A|+j$, for every $b \in
B \setminus A$. Thus
\[
B\setminus A \subseteq \{d-|A|+j, \dots, d-|B|+j+1\}.
\]
Hence
\begin{align*}
|B\setminus A| &\leq \{d-|A|+j, \dots, d-|B|+j+1\}|\\
&= (d-|A|+j)-(d-|B|+j) \\
&= |B|-|A|,
\end{align*}
which implies that $A \subset B$. But $d-|A|+j \notin B$,
otherwise $ \sigma \subset B$. Therefore,
\[
B \subseteq A \cup \{d-|A|+j-1, \dots, d-|B|+j+1\}.
\]
Hence
\begin{align*}
|B| &\leq | A \cup \{d-|A|+j, \dots, d-|B|+j+1\}| \\
&= |A| + (d-|A|+j-1)-(d-|B|+j) \\
&= |B|-1,
\end{align*}
a contradiction.

\noindent {\bf Case 2}. There exists $b \in B\setminus A$ such
that $b > d-|A|+j$, say $b = d-|A|+j +h$ for some $h \geq 1$.
Since $\sigma$ is a non-face and $b > d-|A|+j$, the set $A \cup
\{b\}$ must contain a minimal non-face. Let
\begin{align*}
A^> &:= \{ a \in A \, : \, a > b \} = \{ a_1 > \cdots > a_t\}, \mbox{ and }\\
A^< &:= \{ a \in A \, : \, a < b \} = \{ b_1 > \cdots > b_{|A|-t}
\},
\end{align*}
if $A^> = \emptyset$, let $t:=0$. If $A^< \neq \emptyset$, then
for $1 \leq i \leq |A|-t$,
\[
b_i  = d-|A|+j+s_i,
\]
where $s_i \geq |A|-t-i+1$. Thus
\[
h > s_1 > \cdots > s_{|A|-t} \geq 1.
\]
In the following, we will discuss all possible cases.

\noindent 1. Suppose $A^> \cup \{b\}$ is a minimal non-face of
$\Gamma$. We need to show that $A^> \cup \{b\} \in \Sigma_s$ for
some $s > j$. So it is enough to show that $b = d-|A^>|+s =
d-t+s$, for some $s >j$.  But
\begin{align*}
b > b_1 &= d-|A|+j+s_1 \\
&\geq d-|A|+j + |A| -t \\
&= d-|A^>|+j.
\end{align*}
This implies that $b = d- |A^>|+s$ for some $s > j$. Therefore,
$A^> \cup \{b\} \in \Sigma_s$ and $s >j$.

\vspace{0.1cm} \noindent 2. Suppose  $A^> \cup \{b\} \cup
\{b_1,\dots, b_i\}$ is a minimal non-face of $\Gamma$. As above,
\begin{align*}
b_i &= d-|A|+j+s_i \\
&\geq d-|A|+j + |A|-t-i+1 \\
&= d-(t+i)+j+1.
\end{align*}
Thus $b_i = d-|A^> \cup \{b\} \cup \{b_1,\dots, b_{i-1}\}|+ s$ for
some $s >j$. Therefore,
\[
(A^> \cup \{b\} \cup \{b_1,\dots, b_{i-1}\})\cup\{b_i\} \in
\Sigma_s.
\]
\end{proof}

In the following we construct a generic  deformation of
$I_{\Gamma}$. For $1\leq j \leq n-d$, let
\[
M_j := \langle \prod_{i \in \sigma} x_i^{n-d-j+1}  \, : \, \sigma
\in \Sigma_j \rangle,
\]
a monomial ideal in $R:= k[x_1,\dots,x_n]$. Let
\[
M := M_1 + \cdots + M_{n-d}.
\]
The following corollary is straightforward from Proposition
\ref{thm3:shifted}.
\begin{corollary}\label{coro:shifted}
For any $ 1 \leq j < n-d$, let $m_1, m_2 \in M_j$ be two distinct
generators of $M_j$. Then there exists $ s > j$ and a minimal
generator $ m_3 \in M_s$ such that $m_3$ divides  $\lcm(m_1,
m_2)$. Thus $M$ is a generic  deformation of $I_{\Gamma}$.
\end{corollary}
The following theorem shows that $M$ is Cohen-Macaulay. Before we
prove that, we verify it for our Example.

\vspace{.2cm} \noindent {\it Example \ref{ex1:shifted}} ({\it
cont.}). The above construction yield the following monomial
ideals:
\begin{align*}
M_1 &=  \langle x_2^3x_3^3x_4^3 , x_3^3x_5^3, x_3^3x_6^3\rangle, \\
M_2 &=  \langle x_4^2x_5^2, x_4^2x_6^2\rangle, \\
M_3 &=  \langle x_5x_6 \rangle.
\end{align*}
It is easy to check that the monomial ideal
\[
M = \langle x_2^3x_3^3x_4^3 , x_3^3x_5^3, x_3^3x_6^3, x_4^2x_5^2,
x_4^2x_6^2, x_5x_6 \rangle
\]
is generic and $\rad(M) = I_{\Gamma}$. Moreover, by Corollary
\ref{coro:shifted}, $\dim(\dm) \leq 2$. To show that $M$ is
Cohen-Macaulay, it is enough to show that $\dim(\dm) = 2$, by
Theorem \ref{thethm}. But it is clear that $\lcm\{ x_5x_6,
x_4^2x_5^2, x_3^3x_5^3\}$ is uniquely attained and hence $\{
x_5x_6, x_4^2x_5^2, x_3^3x_5^3\}$ corresponds to a facet of $\dm$.
Therefore  $\dim(\dm) = 2$ and hence $M$ is Cohen-Macaulay.
\begin{theorem}\label{thm5:shifted}
The ideal $M$, constructed above, is a generic  Cohen-Macaulay
integral deformation of  $I_{\Gamma}$.
\end{theorem}
\begin{proof}
By Corollary \ref{coro:shifted}, $M$ is a generic integral
deformation of $I_{\Gamma}$. So it remains to show that $M$ is
Cohen-Macaulay. By Theorem \ref{thethm}, it is enough to show that
$\dim(\dm) = n-d-1$. By Corollary  \ref{coro:shifted}, $\dim(\dm)
\leq n-d-1$. Thus to show that $\dim(\dm) = n-d-1$, it is enough
to find a face $F$ of $\dm$ such that $\dim(F) = n-d-1$.

For $ 1 \leq j \leq n-d$, let
\[
\sigma_{j} : = \{d+j, d+j-1,\dots,i_j\},
\]
where $i_j$ is small enough such that $\sigma_j \in
\Sigma_{\Gamma}$, i.e., $\sigma_j$ is a minimal non-face of
$\Gamma$. Such an $i_j$ exists since $|[d+j]| > d$ for $j >1$ and
$\dim(\Gamma)=d-1$.
 First we will show that $\sigma_j \in \Sigma_j$.
Let $A : = \sigma_j \setminus \{i_j\}$. It is clear that $|A| =
(d+j)-i_j$. Thus $i_j = d -|A| +j$. Therefore $A \subseteq
\overline{[d-|A|+j]}$ and hence $\sigma = A \cup  \{i_j\} \in
\sigma_j$.

Next we will show that $F:=\{m_{\sigma_1}, \dots,
m_{\sigma_{n-d}}\}$ is a face of $\dm$. Now
\begin{align*}
\lcm(F) &= \lcm(m_{\sigma_1}, \dots, m_{\sigma_{n-d}})\\
&= x_nx_{n-1}^2 \cdots x_{d+1}^{n-d}x_{d}^{n-d} \cdots
x_{i_1}^{n-d}.
\end{align*}
It is clear that, for $1 \leq j \leq n-d$, we have $\lcm(F) \neq
\lcm(F\setminus \{m_{\sigma_j}\})$. So we only need to show that
$\lcm(F) \neq  \lcm(F \cup \{m_{\sigma}\})$ for any $\sigma \in
\Sigma_{\Gamma}\setminus F$. Let $\sigma \in
\Sigma_{\Gamma}\setminus F$, say $\sigma \in \Sigma_j$ for some $1
\leq j \leq n-d$. Thus $\sigma = A \cup \{d-|A|+j\}$, for some $A
\in \Gamma$ and $A \subseteq \overline{[d-|A|+j]}$. Since $\sigma
\notin F$, there exists $a \in A$ such that $a > d+j$, otherwise
$\sigma = \sigma_j$. Suppose $a = d+i$ for some $i >j$. Thus $a
\in \sigma_i$ and hence $\deg_{x_a}(\lcm(F)) = n-(d+i)+1$. But
$\deg_{x_a}(m_{\sigma})= n-(d+j)+1$. Therefore, $\lcm(F) \leq
\lcm(F \cup \{m_{\sigma}\})$. Hence $F$ is a face of $\dm$ of
dimension $n-d-1$. Therefore, $\dim(\dm) = n-d-1$. By Theorem
\ref{thethm}, $M$ is Cohen-Macaulay.
\end{proof}

\section{Tree Complexes: Shellable Clique Complexes of Chordal Graphs}
\hspace{\parindent} In this section, we study the class of
shellable clique complexes of chordal graphs,  which we call tree
complexes because they generalize tree graphs. The main result of
this section is  that tree complexes have generic Cohen-Macaulay
deformations. Moreover,  we  give an exact formula to compute the
$f$-vector of tree complexes.

Throughout this section, let $G$ be a simple graph (undirected
graph with no loops and no multiple edges) on the vertex set
$[n]$.
\begin{definition}\label{def1:flag}
Let $G$ be as above.
\begin{enumerate}
\item A subset $\sigma \subseteq [n]$ is a {\it clique} if the
induced subgraph $G_\sigma$ is complete. \item The \emph{clique
complex} of $G$ is
\[
K(G) := \{\sigma \subseteq [n] \, : \, \sigma \mbox{ is a clique
of } G\}.
\]
\end{enumerate}
\end{definition}
\begin{example}\label{ex1:flag}
A graph $G$ and its clique complex $K(G)$ are in Figure
\ref{fig1:flag}. Notice that $G$ is the 1-skeleton of $K(G)$.
\begin{figure}[!htp]
  \centering
  \includegraphics[width=10cm]{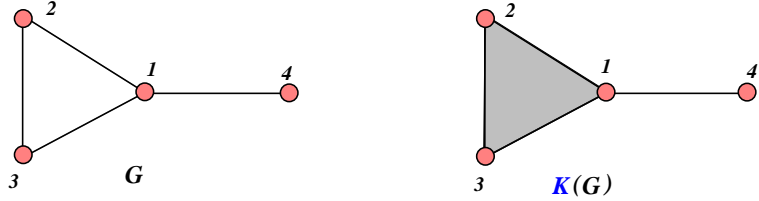}
  \caption{A graph $G$ and its clique complex $K(G)$.}
  \label{fig1:flag}
\end{figure}
\end{example}
The proof of the following lemma is straightforward.
\begin{lemma}\label{lem1:flag}
Every minimal non-face of $K(G)$ has precisely two elements.
Moreover, if $\Delta$ is a simplicial complex such that every
minimal non-face has exactly two elements, then $\Delta=K(G)$,
where $G$ is the 1-skeleton of $\Delta$.
\end{lemma}
Clique complexes are also called \emph{flag} complexes, see
\cite{St1}.

\vspace{.2cm}
\begin{definition}\label{def3:flag}
A simplicial complex $\Delta$ is a {\it tree} complex if there
exists a shelling $F_1, \dots, F_t$ of $\Delta$ such that, for $1
< i \leq t$, there exists $j_i \in F_i$ such that $j_i \notin F_1
\cup \cdots \cup F_{i-1}$. In this case, we say that $F_i$
introduces the vertex $j_i$.
\end{definition}
It is easy to see that any tree graph is a 1-dimensional tree
complex. Hence tree complexes are a natural generalization of
tree graphs.
\begin{example}\label{ex-2:tree}
The simplicial complex $\Gamma$ in Figure \ref{fish1} is a
2-dimensional tree complex.
\begin{figure}[!htp]
  \centering
  \includegraphics[height=12cm]{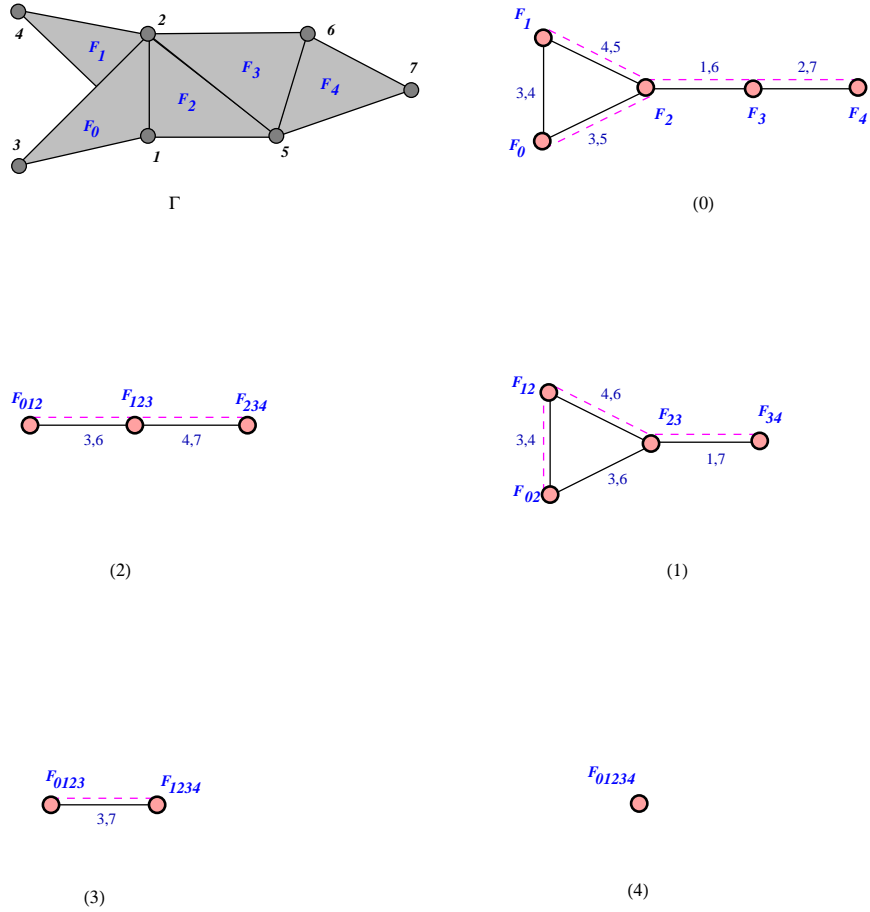}
  \caption{The tree complex from Example \ref{ex-2:tree}.}
  \label{fish1}
\end{figure}
\end{example}
One of the most interesting problems in combinatorics is
characterizing the  $f$-vectors of flag complexes \cite{BB, St1,
K1}. For tree complexes,
 the next theorem gives
an exact formula for the $f$-vector of these complexes.
\begin{theorem}\label{thm3:flag}
Let $\Delta$ be a $(d-1)$-dimensional tree complex on $[n]$. Then
the number of facets of $\Delta$ is $f_{d-1}(\Delta) = n-d+1$, and
for $ 1 < i \leq d-1$,
\begin{align*}
f_{i-1}(\Delta) &:= {d\choose i} + (n-d)\lbrace{d\choose i} - {{d-1}\choose i}\rbrace\\
&= (n-d+1){d \choose i } - (n-d) {{d-1}\choose i}.
\end{align*}
\end{theorem}
\begin{proof}
Let $F_1, \dots, F_t$ be a shelling of $\Delta$. By Theorem
\ref{thm2:flag}, every vertex $d \leq i \leq n$ introduces a
facet. Thus $f_{d-1}(\Delta) = n-d+1$. Now for $1 < i \leq d-1$,
it is clear that $F_1$ has $d \choose i$ $(i-1)$-faces. Moreover,
every facet $F_j$, where $j >1$ has $d \choose i$ $(i-1)$-faces
minus the the number of  $(i-1)$-faces that $F_j$ shares with
$F_1, \dots, F_{j-1}$. But $F_j$ shares only
 $d-1$ vertices with  $F_1, \dots, F_{j-1}$. Hence the number of
$(i-1)$-faces it shares is ${d-1} \choose i$ which implies the
equality above.
\end{proof}
\begin{definition}\label{def2:flag}
A graph $G$ is {\it chordal} if every cycle of length four or more
has a chord (an edge linking two non-adjacent nodes in the cycle).
\end{definition}

\begin{example}
The graph $G$ in Figure \ref{fig1:flag} is chordal.
\end{example}
The following theorem gives many useful characterizations of
chordal graphs.
\begin{theorem}[{\cite[Theorem 7, p 112]{P}}]\label{thm1:flag}
The following statements are equivalent.
\begin{enumerate}
\item $G$ is a chordal graph. \item All vertices of $G$ can be
deleted by arranging them in separate piles, one for each maximal
clique, and then repeatedly applying the following two operations:
\begin{itemize}
\item Delete a vertex that occurs in only one pile. \item Delete a
pile if all its vertices appear in another pile.
\end{itemize}
\item There is a spanning tree (called a clique tree) $T$ of the
facet graph of $K(G)$ such that for every vertex $i$ of $G$, if we
remove from $T$ all cliques not containing $i$, the remaining
subtree stays connected. In other words, any two cliques
containing $i$ are either adjacent in $T$ or connected by a path
made entirely of cliques that contain $i$.
\end{enumerate}
\end{theorem}
Suppose $\Delta$ is a $(d-1)$-dimensional tree complex on $[n]$.
Without loss of generality (after relabeling the vertices), we may
assume $F_1 : =[d]$ and, for $1 < i \leq t$, $F_i$ introduces the
vertex $d+i-1$, i.e., $F_i : = A \cup \{d+i-1\}$ for some $A
\subseteq [d+i-2]$ and $|A| = d-1$. The following theorem shows
that the class of tree complexes is equal to the class of
shellable clique complexes of chordal graphs.
\begin{theorem}\label{thm2:flag}
A shellable simplicial complex $\Delta$ is a tree complex if and
only if $\Delta = K(G)$, the clique complex of a  chordal graph
$G$.
\end{theorem}
\begin{proof}
Let $F_1, \dots, F_t$ be a  shelling of $\Delta$. First suppose
that $\Delta = K(G)$, for some chordal graph $G$, and assume that
$\Delta$ is not a tree complex. Let $s$ be the smallest index such
that $F_s$ does not introduce a new vertex. There exists a simple
cycle in the facet graph of $K(G)$ that contains $F_s$, without
loss of generality, $F_1,\dots,F_s$ is a simple cycle of minimal
length. Thus every vertex is contained in more than one clique.
This is a contradiction to (2) of  Theorem \ref{thm1:flag}.
Therefore, for $1 < i \leq t$, there exists $i_j \in F_i$ such
that $i_j \notin F_1 \cup \cdots \cup F_{i_1}$ and hence $\Delta$
is a tree complex.

Now suppose that $\Delta$ is a tree complex. Let $G$ be the
1-skeleton of $\Delta$ ($G$ is the subcomplex of simplexes of
dimension less than or equal to 1). It is clear that $G$ is a
graph on the set of vertices $[n]$. Indeed $G$ is a chordal graph:
For every facet $F_i$, let $P_i$ be the pile that contains the
vertices of $F_i$. Since, for $1 < i \leq t$, $F_i$ introduces the
vertex $j_i$, it is clear that one can empty all the piles by
deleting each vertex that is in only one pile and then deleting
any pile where its vertices are in another pile. Thus $G$ is
chordal. So we only need to show that $\Delta = K(G)$. It is clear
that $\Delta$ is a subcomplex of $K(G)$. Suppose $C$ is a clique
of $G$ and $h \in C$ is the maximal element. Thus the vertex $h$
is introduced by the facet $F_{i_h}$, for some $1 \leq i_h \leq
t$. Thus every vertex connected to $h$ and less than $h$ is in
$F_{i_h}$ which implies that $C \subset F_{i_h}$. Hence $\Delta =
K(G)$.
\end{proof}
Let $FG(\Delta)$ be the facet graph of $\Delta$. By Theorem
\ref{thm3:flag},  $FG(\Delta)$ has $n-d+1$ facets. Let
$T_{\Delta}$ be a clique tree of $FG(\Delta)$, it has $n-d$ edges
and each edge in $T_{\Delta}$  corresponds to a missing edge in
$\Delta$:

If $F_i$ and $F_j$  are vertices connected by an edge $e$ in
$T_{\Delta}$, then $F_i$ and $F_j$ are facets of $\Delta$, i.e.,
$|F_i|=|F_j|= d$ and $|F_i \cap F_j|= d-1$. Thus, since $\Delta$
is a clique complex,  there is a unique minimal non-face in $F_i
\cup F_j$, say $\{p_e,q_e\}$.  Let
\[
\Sigma_0 := \{ \{p_e,q_e\} \, : \, e \in T_{\Delta} \}.
\]
Now for the edge $e$, let
\[
F_{ij} := F_i \cup F_j.
\]
Let $\Delta_1$ be the $d$-dimensional simplicial complex on $[n]$
with the set of facets $\{ F_{ij} \, : \,  e \in T_{\Delta}\}$. It
is easy to see that $\Delta_1$ is a tree complex. If $\Delta_1$ is
not a simplex, let $T_{\Delta_1}$ be a clique tree of
$FG(\Delta_1)$. Repeat the steps above for this case ($\Delta_t$
is a simplex only when $t = n-d$). Hence, in general, for $1 \leq
t < n-d$, let
\[
\Sigma_t := \{\{p_e,q_e\} \, : \, e \in T_{\Delta_t}\}.
\]
This process gives a partition of  $\Sigma_{\Gamma}$. Before we
give a proof of this fact, we will verify it for our example.

\vspace{.2cm} \noindent \emph{Example} \ref{ex-2:tree} ({\it
cont.}). It is easy to check that
\[
\Sigma_{\Gamma} =
\{\{3,5\},\{4,5\},\{1,6\},\{2,7\},\{3,4\},\{4,6\},\{1,7\},\{3,6\},\{4,7\},\{3,7\}\}.
\]
Moreover, the graph in Figure \ref{fish1}~($i$) is $FG(\Gamma_i)$
and the subgraph of dashed edges is $T_{\Gamma_i}$.
This yields the following partition of $\Sigma_{\Gamma}$.
\begin{align*}
\Sigma_0 &= \{\{3,5\},\{4,5\},\{1,6\},\{2,7\}\} \\
\Sigma_1 &= \{\{3,4\},\{4,6\},\{1,7\},\{3,6\},\{4,7\}\}\\
\Sigma_2 &= \{\{3,6\},\{4,7\}\}\\
\Sigma_3 &= \{\{3,7\}\}.
\end{align*}
The proof of the following proposition is straightforward.
\begin{proposition}
Let $\Delta$ be a tree complex and $\Sigma_t$ be as above.
\begin{enumerate}
\item For $0 \leq t < n-d$,  we have $\Sigma_t \neq \emptyset$,
\item For $0 \leq t_1 < t_2 < n-d$, $\Sigma_{t_1} \cap
\Sigma_{t_2} = \emptyset$, and \item  $\Sigma_{\Delta} =
\cup_{t=0}^{n-d-1} \Sigma_t$.
\end{enumerate}
\end{proposition}
The following proposition
 will help us to construct a generic deformation from $I_{\Delta}$.
\begin{proposition}\label{prop-3:flag}
Let $\Delta$, $\Sigma_t$ be as above. Let $e_1, e_2 \in \Sigma_t$
and $e_1  \neq e_2$. There exists $ s > t$ and $ e_3 \in \Sigma_s$
such that $e_3 \subset e_1 \cup e_2$.
\end{proposition}
\begin{proof}
Without loss of generality, let
\begin{align*}
e_1 = \{p_1,q_1\}, \qquad e_2 = \{p_2,q_2\}
\end{align*}
be two edges of the clique tree $T_{\Delta_t}$ where the edge
$e_1$ (resp. $e_2$) has as vertices the $(d+t-1)$-simplexes
(facets) $\tau_1,\tau_2$ (resp. $\gamma_1,\gamma_2$) of
$\Delta_t$, say
\begin{align*}
\tau_1 = A \cup \{p_1\}, \qquad \tau_2 = A \cup \{q_1\},\\
\gamma_1 = B \cup \{p_2\}, \qquad \gamma_2 = B \cup \{q_2\}.
\end{align*}
{\bf Case 1}. $p := p_1 = p_2$ ($q_1 \neq q_2$). If
$\tau_1=\gamma_1$, then $A = B$. Thus the edge $\{q_1,q_2\}$ is in
the facet graph $FG(\Delta_t)$ but not in $T_{\Delta_t}$,
otherwise $\tau_1, \tau_2, \gamma_2$ form a cycle in
$T_{\Delta_t}$. So assume that $\tau_1 \neq \gamma_1$. There
exists a path $\tau_1=G_1, \dots,G_r = \gamma_1$ such that $p \in
G_i$ for $1 \leq i \leq r$. This implies that $\tau_2,G_1,
\dots,G_r,\gamma_2$ is the only
 path between $\tau_2$ and $\gamma_2$ in  $T_{\Delta_t}$. Hence
$\{q_1,q_2\}$ is not an edge in  $T_{\Delta_t}$, otherwise we have
a cycle in  $T_{\Delta_t}$. In this case, let $e_3 = \{q_1,q_2\}$.

{\bf Case 2}. $p_i \neq q_j$ for $1 \leq i, j \leq 2$. Suppose for
$1 \leq i, j \leq 2$, the edge $\{p_i,q_j\} \in \Delta_t$. Thus
$p_1,q_1,p_2,q_2$ form a 4-cycle in the chordal graph $G$, hence
either the edge $\{p_1,q_1\}$ or $\{p_2,q_2\}$ is in $G$, a
contradiction. Without loss of generality, suppose the 1-simplex
$\{p_1,p_2\} \notin \Delta_t$. If $\{p_1,p_2\} \notin
T_{\Delta_t}$, then we are done. Suppose $\{p_1,p_2\}$ is an edge
in $T_{\Delta_t}$ with the vertices $F$ and $G$ such that $p_1 \in
F$ and $p_2 \in G$. Since $p_1 \in F \cap \tau_1$, there exists a
path between $\tau_1$ and $F$ made entirely from facets containing
$p_1$. Similarly, there exists a path between $G$ and $\gamma_1$
made entirely from facets containing $p_2$. Thus there is a path
between $\tau_2$ and $\gamma_2$. It is clear that if $\{q_1,q_2\}
\in \Delta_t$ or  $\{p_1,p_2\}$ is an edge in $T_{\Delta_t}$, then
we have a cycle in $T_{\Delta_t}$. Once again, in this case, let
$e_3 = \{q_1,q_2\}$.
\end{proof}
In the following we construct a generic integral deformation of
$I_{\Delta}$. For $0 \leq t < n-d$, and $\{p_e,q_e\} \in
\Sigma_t$, let
\[
M_t := \langle x_{p_e}^{n-d-t}x_{q_e}^{n-d-t} \, : \,  \{p_e,q_e\}
\in \Sigma_t\rangle.
\]
It is clear that $M_t \neq 0$, otherwise there are no missing
faces and hence we have a simplex which is not the case unless $t
= n-d$.

Let $M = M_0 + M_1 +\cdots + M_{d-t-1}$. It is clear that $\rad(M)
= I_{\Delta}$. Moreover, in the next Lemma, we show that $M$ is
generic. Before that we continue with our example.

\vspace{.2cm} \noindent \emph{Example} \ref{ex-2:tree} ({\it
cont.}). For each $i$; $0 \leq i \leq 4$, the corresponding ideals
to the partition above are:
\begin{align*}
M_0 &= \langle x_3^4x_5^4,x_4^4x_5^4,x_1^4x_6^4,x_2^4x_7^4\rangle, \\
M_1 &= \langle x_3^3x_4^3,x_4^3x_6^4, x_1^3x_7^3\rangle, \\
M_2 &= \langle x_3^2x_6^2,x_4^2x_7^2\rangle, \\
M_3 &= \langle x_3x_7\rangle.
\end{align*}
Since $\Gamma_4$ is a simplex,  the facet graph $FG(\Gamma_4)$ is
a vertex with no edges and hence $M_4 = 0$.

Let $M = M_0 + \cdots +M_3$. It is easy to see that $M$ is a
generic monomial ideal and $V(M)= \Gamma$.

The following lemma follows from Corollary \ref{prop-3:flag}.
\begin{lemma}\label{lemma-flag}
For any $ 0 \leq j < n-d$, let $m_1, m_2 \in M_j$ be two distinct
generators of $M_j$. Then there exists $ s > j$ and a minimal
generator $ m_3 \in M_s$ such that $m_3$ divides  $\lcm(m_1,
m_2)$. Thus $M$ is a generic integral deformation of $I_{\Delta}$.
\end{lemma}
The following theorem shows that $M$ is Cohen-Macaulay. Before we
prove that, we verify it for our example.

\vspace{.2cm} \noindent \emph{Example} \ref{ex-2:tree} ({\it
cont.}). It is easy to see that
$\{x_3x_7,x_3^2x_6^2,x_3^3x_4^3,x_3^4x_5^4\}$ is a face of $\dm$
and hence $\dim(\dm)= 3=7-3-1$. Thus, by Theorem \ref{thethm}, $M$
is Cohen-Macaulay.

\begin{theorem}\label{thm5:flag}
The ideal $M$ is a generic Cohen-Macaulay integral deformation of
$I_{\Delta}$.
\end{theorem}
\begin{proof}
It remains to show that $M$ is Cohen-Macaulay. Since
$V(M)=\Delta$, by Theorem \ref{thethm}, we need to show that
$\dim(\dm)= n-d-1$. By Lemma \ref{lemma-flag}, each face of $\dm$
is of dimension $ \leq n-d-1$. To show that  $\dim(\dm)= n-d-1$,
it is enough to find a face of dimension $n-d-1$. For $0 \leq i <
n-d$, there exists a unique minimal non-face in $\Delta_{i}$ which
is a subset of  $[d+i+1]$ and contains $d+i+1$, say
$\{\alpha_i,d+i+1\} \in \Sigma_{i}$. Let $m_i=
x_{\al_i}^{n-d-i}x_{d+i+1}^{n-d-i}$ be the corresponding monomial
in $M_{i}$.

\noindent {\underline {Claim}}. $C := \{m_0, \dots,m_{n-d-1}\}$ is
a facet of $\dm$ and hence $\dim(\dm) = n-d-1$.
\[
\lcm(C):= \lcm\{m_0,\dots,m_{n-d-1}\}= f\cdot
x_{d+1}^{n-d}x_{d+2}^{n-d-1} \cdots x_{d+i}^{n-d-i+1} \cdots
x_{n-1}^2x_n,
\]
where $f$ is a monomial on the variables $x_1,\dots,x_{d}$. It is
clear that $\lcm(C) \neq \lcm(H)$ for any proper subset $H \subset
C$. Suppose there exists $i,j$ such that $j < d+i+1$, $\{j,d+i+1\}
\in \Sigma_{\Delta}$ and $x_j^{a}x_{d+i+1}^{a}$ divides $\lcm(C)$,
i.e., $a \leq n-d-i+2$. This implies that $\{j,d+i+1\}$
corresponds to an edge in $\Delta_{s}$ where $s \geq i$. Since
$\{\alpha_e,d+i+1\}$ is the only non-face in $\Delta_{i}$ that is
a subset of $[d+i+1]$, we must have $s=i$ and
 $j = {\alpha_e}$ and hence $m = m_{i}$.
Therefore, $C \in \dm$ and hence $\dim(\dm) = n-d-1$. By Theorem
\ref{thethm}, $M$ is Cohen-Macaulay.
\end{proof}

\bibliography{paper}
\bibliographystyle{plain}
\end{document}